\documentclass{amsart}

\allowdisplaybreaks
\theoremstyle{plain}
\newtheorem{thm}{Theorem}[section]

\newtheorem{lem}[thm]{Lemma}

\theoremstyle{definition}

\theoremstyle{remark}

\newtheorem*{remark*}{Remark}

\newcommand{\cF}{{\mathcal{F}}}

        \newcommand{\field}[1]{{\mathbb{#1}}}
        \newcommand{\NN}{\field{N}}
        \newcommand{\ZZ}{\field{Z}}
        
        \newcommand{\RR}{\field{R}}

\newcommand{\Dom}{\mbox{\rm Dom}}

\newcommand{\Tr}{\mbox{\rm Tr}}

\title[The periodic magnetic Schr\"odinger operators]{The periodic magnetic Schr\"odinger operators: spectral gaps
and tunneling effect}

\author{Bernard Helffer}
\address{D\'epartement de Math\'ematiques, B\^atiment 425, Universit\'e Paris-Sud,
F91405 Orsay C\'edex, France} \email{Bernard.Helffer@math.u-psud.fr}

\author{Yuri A. Kordyukov}
\address{Institute of Mathematics, Russian Academy of Sciences, 112 Chernyshevsky str., 450077
Ufa, Russia} \email{yurikor@matem.anrb.ru}
\thanks{B.H. was partially supported by the ESF programme SPECT.
Y.K. was partially supported by the Russian Foundation of Basic
Research (grant 06-01-00208).}

\date{}

\begin{document}
\bibliographystyle{plain}

\begin{abstract}
A periodic Schr\"odinger operator on a noncompact Riemannian
manifold $M$ such that $H^1(M, \RR)=0$ endowed with a properly
discontinuous cocompact isometric action of a discrete group is
considered. Under some additional conditions on the magnetic field
existence of an arbitrary large number of gaps in the spectrum of
such an operator in the semiclassical limit is established. The
proofs are based on the study of the tunneling effect in the
corresponding quantum system.
\end{abstract}

\maketitle

\section*{Introduction}
Let $ M$ be a noncompact oriented manifold of dimension $n\geq 2$
equipped with a properly discontinuous action of a finitely
generated, discrete group $\Gamma$ such that $M/\Gamma$ is
compact. Suppose that $H^1(M, \RR) = 0$, i.e. any closed $1$-form
on $M$ is exact. Let $g$ be a $\Gamma$-invariant Riemannian metric
and $\bf B$ a real-valued $\Gamma$-invariant closed 2-form on $M$.
Assume that $\bf B$ is exact and choose a real-valued 1-form $\bf
A$ on $M$ such that $d{\bf A} = \bf B$.

Thus, one has a natural mapping
\[
u\mapsto ih\,du+{\bf A}u
\]
from $C^\infty_c(M)$ to the space $\Omega^1_c(M)$ of smooth,
compactly supported one-forms on $M$. The Riemannian metric allows
to define scalar products in these spaces and consider the adjoint
operator
\[
(ih\,d+{\bf A})^* : \Omega^1_0(M)\to C^\infty_0(M).
\]
A Schr\"odinger operator with magnetic potential $\bf A$ is
defined by the formula
\[
H^h = (ih\,d+{\bf A})^* (ih\,d+{\bf A}).
\]
Here $h>0$ is a semiclassical parameter, which is assumed to be
small.

Choose local coordinates $X=(X_1,\ldots,X_n)$ on $M$. Write the
1-form $\bf A$ in the local coordinates as
\[
{\bf A}= \sum_{j=1}^nA_j(X)\,dX_j,
\]
the matrix of the Riemannian metric $g$ as
\[
g(X)=(g_{jl}(X))_{1\leq j,l\leq n}
\]
and its inverse as
\[
g(X)^{-1}=(g^{jl}(X))_{1\leq j,l\leq n}.
\]
Denote $|g(X)|=\det(g(X))$. Then the magnetic field $\bf B$ is
given by the following formula
\[
{\bf B}=\sum_{j<k}B_{jk}\,dX_j\wedge dX_k, \quad
B_{jk}=\frac{\partial A_k}{\partial X_j}-\frac{\partial
A_j}{\partial X_k}.
\]
Moreover, the operator $H^h$ has the form
\[
H^h=\frac{1}{\sqrt{|g(X)|}}\sum_{1\leq j,l\leq n}\left(i h
\frac{\partial}{\partial X_j}+A_j(X)\right)\left[\sqrt{|g(X)|}
g^{jl}(X) \left(i h \frac{\partial}{\partial
X_l}+A_l(X)\right)\right].
\]

For any $x\in M$, denote by $B(x)$ the anti-symmetric linear
operator on the tangent space $T_x{ M}$ associated with the 2-form
$\bf B$:
\[
g_x(B(x)u,v)={\bf B}_x(u,v),\quad u,v\in T_x{ M}.
\]
Recall that the intensity of the magnetic field is defined as
\[
{\Tr}^+ (B(x))=\sum_{\substack{\lambda_j(x)>0\\ i\lambda_j(x)\in
\sigma(B(x)) }}\lambda_j(x)=\frac{1}{2}\Tr([B^*(x)\cdot
B(x)]^{1/2}).
\]
It turns out that in many problems the function $x\mapsto h\cdot
{\Tr}^+ (B(x))$ can be considered as a magnetic potential, that is,
as a magnetic analogue of the electric potential $V$ in a
Schr\"odinger operator $-h^2\Delta+V$.

In this paper we will always assume that the magnetic field has a
periodic set of compact potential wells. More precisely, put
\[
b_0=\min \{{\Tr}^+ (B(x))\, :\, x\in { M}\}
\]
and assume that there exist a (connected) fundamental domain $\cF$
and a constant $\epsilon_0>0$ such that
\begin{equation}\label{e:tr1}
  {\Tr}^+ (B(x)) \geq b_0+\epsilon_0, \quad x\in \partial\cF.
\end{equation}
For any $\epsilon_1 \leq \epsilon_0$, put
\[
U_{\epsilon_1} = \{x\in \cF\,:\, {\Tr}^+ (B(x)) < b_0+
\epsilon_1\}.
\]
Thus $U_{\epsilon_1}$ is an open subset of $\cF$ such that
$U_{\epsilon_1}\cap \partial\cF=\emptyset$ and, for $\epsilon_1 <
\epsilon_0$, $\overline{U_{\epsilon_1}}$ is compact and included in
the interior of $\cF$. Any connected component of $U_{\epsilon_1}$
with $\epsilon_1 < \epsilon_0$ and also any its translation under
the action of an element of $\Gamma$ can be understood as a magnetic
well. These magnetic wells are separated by potential barriers,
which are getting higher and higher when $h\to 0$ (in the
semiclassical limit).

For any linear operator $T$ in a Hilbert space, we will denote by
$\sigma(T)$ its spectrum. By a gap in the spectrum of a
self-adjoint operator $T$ we will mean any connected component of
the complement of $\sigma(T)$ in $\RR$, that is, any maximal
interval $(a,b)$ such that
\[
(a,b)\cap \sigma(T) = \emptyset\,.
\]

We will consider the magnetic Schr\"odinger operator $H^h$ as an
unbounded self-adjoint operator in the Hilbert space $L^2(M)$. The
main object of our investigation is the gaps in the spectrum of this
operator.

The problem of existence of gaps in the spectra of second order
periodic differential operators has been extensively studied
recently. Some related results on spectral gaps for periodic
magnetic Schr\"odinger operators can be found for example in
\cite{BDP,gaps,HS88,HSLNP345,HempelHerbst95,HempelPost02,HerbstNakamura,KMS,Ko04,bedlewo,MS,Nakamura95}
(see also the references therein).

In this paper, we will be interested in spectral gaps located below
the top of potential barriers, that is, on the interval $[0,
h(b_0+\epsilon_0)]$. For such energy levels of a quantum particle,
the important role is played by the tunneling effect, that is, by
the possibility for a particle with such an energy to pass through a
potential barrier. We will demonstrate how the investigation of the
tunneling effect allows to obtain some information on gaps in the
spectrum of the Schr\"odinger operator in the given energy interval.

The results of the paper were partially obtained during our stay
in the Isaac Newton Institute of Mathematical Sciences, Cambridge,
to whom we express our gratitude for hospitality and support.

\section{Tunneling effect and localization of the spectrum}\label{s:tunnel}
In this Section, we briefly recall the key result on localization
of the spectrum of the magnetic Schr\"odinger operator $H^h$
obtained in \cite{gaps}, which follows from the semiclassical
analysis of the tunneling effect for the corresponding quantum
system.

For any domain $W$ in $ M$, denote by $H^h_W$ the unbounded
self-adjoint operator in the Hilbert space $L^2(W)$ defined by the
operator $H^h$ in $\overline{W}$ with the Dirichlet boundary
conditions. The operator $H^h_W$ is generated by the quadratic form
\[
u\mapsto q^h_W [u] : = \int_W |(ih\,d+{\bf A})u|^2\,dx
\]
with the domain
\[
\Dom (q^h_W) = \{ u\in L^2(W) : (ih\,d+{\bf A})u \in L^2\Omega^1(W),
u\left|_{\partial W}\right.=0 \},
\]
where $L^2\Omega^1(W)$ denotes the Hilbert space of $L^2$
differential $1$-forms on $W$, $dx$ is the Riemannian volume form
on $M$.

Assume now that the operator $H^h$ satisfies the condition
(\ref{e:tr1}). Fix $\epsilon_1>0$ and $\epsilon_2>0$ such that
$\epsilon_1 < \epsilon_2 < \epsilon_0$, and consider the operator
$H^h_D$ associated with the domain $D=\overline{U_{\epsilon_2}}$.
The operator $H^h_D$ has discrete spectrum.

The investigation of the tunneling effect allows to prove that the
spectrum of $H^h$ on the interval $[0, h(b_0+\epsilon_1)]$ is
localized in an exponentially small neighborhood of the spectrum
of $H^h_D$.

\begin{thm}[\cite{gaps}]\label{t:D}
Under the assumption (\ref{e:tr1}), for any $\epsilon_1 <
\epsilon_2< \epsilon_0$, there exist $C, c, h_0>0$ such that, for
any $h\in (0,h_0]$
\begin{gather*}
\sigma(H^h)\cap [0, h(b_0+\epsilon_1)] \subset \{\lambda\in [0,
h(b_0+\epsilon_2)] : {\rm dist}(\lambda, \sigma(H^h_D))<
Ce^{-c/\sqrt{h}}\},\\ \sigma(H^h_D)\cap [0, h(b_0+\epsilon_1)]
\subset \{\lambda\in [0, h(b_0+\epsilon_2)] : {\rm dist}(\lambda,
\sigma(H^h))< Ce^{-c/\sqrt{h}}\}.
\end{gather*}
\end{thm}

A slightly weaker version of this theorem (which uses the largest
absolute value of the eigenvalues of $B(x)$ instead of ${\Tr}^+
(B(x))$) was proved by Nakamura in \cite{Nakamura95}.

The proof of Theorem~\ref{t:D} uses the approach to the study of the
tunneling effect in multi-well problems developed by Helffer and
Sj\"ostrand for Schr\"odinger operators with electric potentials
(see for instance \cite{HSI,HSII}) and extended to magnetic
Schr\"odinger operators in \cite{HS87,HM}. Since the operator $H^h$
is not with compact resolvent, we do work not with individual
eigenvalues, but with resolvents and use the strategy developed in
\cite{HSII,HS88,Di-Sj} for the case of electric potential and in
\cite{Frank} for the case of magnetic field.

The idea of the proof is to construct an approximate resolvent
$R^h(z)$ of the operator $H^h$ for any $z$, which is not
exponentially close to the spectrum of $H^h_D$, starting from the
resolvent of $H^h_D$ and the resolvent of the Dirichlet
realization of $H^h$ in the complement to the wells. The proof of
the fact that the error of the approximation is exponentially
small is based on Agmon-type weighted estimates (cf. \cite{Ag} and
their semi-classical versions in \cite{HSI} for the case of
Schr\"odinger operators and \cite{HM} for the case of magnetic
Schr\"odinger operators).

\section{Quasimodes and spectral gaps}
Theorem~\ref{t:D} reduces the investigation of gaps in the
spectrum of the operator $H^h$ to the study of the eigenvalue
distribution for the operator $H^h_D$. Actually, it turns out
that, in order to show the existence of arbitrary large number of
gaps in the spectrum of $H^h$ on some interval, it suffices to
construct arbitrarily long sequences of approximate eigenvalues of
$H^h_D$ on this interval located far enough from each other. This
observation is formulated more precisely in the following theorem.

\begin{thm} \label{t:abstract}
Let $N\geq 1$. Suppose that there is a subset
$\mu_0^h<\mu_1^h<\ldots <\mu_{N}^h$ of an interval $I(h)\subset
[0, h(b_0+\epsilon_1))$ such that
\begin{enumerate}
  \item There exist constants $c>0$ and $M \geq 1$ such that
\begin{equation}
\begin{split}\label{e:mu}
& \mu_{j}^h-\mu_{j-1}^h>ch^M, \quad j=1,\ldots,N,\\
& {\rm dist}(\mu_0^h,\partial I(h))>ch^M,\quad {\rm
dist}(\mu_N^h,\partial I(h))>ch^M,
\end{split}
\end{equation}
for any $h>0$ small enough;
  \item Each $\mu_j^h, j=0,1,\ldots,N,$ is an approximate eigenvalue of
the operator $H^h_D$: for some $v_j^h\in C^\infty_c(D)$ we have
\begin{equation}\label{e:eigen}
\|H^h_Dv_j^h-\mu^h_jv_j^h\|=\alpha_j(h)\|v_j^h\|,
\end{equation}
where $\alpha_j(h)=o(h^M)$ as $h\to 0$.
\end{enumerate}
Then the spectrum of $H^h$ on the interval $I(h)$ has at least $N$
gaps for any sufficiently small $h>0$.
\end{thm}

\begin{proof}
Recall the following well-known estimate, which holds for any
self-adjoint operator $A$ in a Hilbert space:
\[
\|(A-\lambda I)^{-1}\|=1/d(\lambda,\sigma(A)), \quad \lambda\not\in
\sigma(A).
\]
By this fact and (\ref{e:eigen}), it follows that, for any
$j=0,1,\ldots,N$, there exists $\lambda_j^h \in \sigma(H^h_D)\cap
I(h)$ such that
\begin{equation}\label{e:mu1}
\lambda_j^h-\mu_j^h=o(h^M), \quad h\to 0.
\end{equation}
By (\ref{e:mu}) and (\ref{e:mu1}), we have
\[
\lambda_{j}^h-\lambda_{j-1}^h>c h^M, \quad j=1,\ldots,N,
\]
for any $h>0$ small enough.

Recall also a rough estimate for the number $N_h(\alpha,\beta)$ of
ei\-gen\-va\-lu\-es of $H_D^h$ on an arbitrary interval
$(h\alpha,h\beta)$ (see, for instance, \cite[Lemma 4.2]{HM}): for
some $C$ and $h_0$
\begin{equation}\label{l:Nh}
  N_h(\alpha,\beta) \leq Ch^{-n}, \quad \forall h\in (0,h_0]\;.
\end{equation}

\begin{lem}\label{l:s}
Let $M>0$ and $c>0$. There exist $C>0$ and $h_1>0$ such that, if
$\alpha^h$ and $\beta^h$ are two points in the spectrum of $H^h$ on
the interval $I(h)$ with $\beta^h-\alpha^h>ch^M$, then, for any
$h\in (0,h_1]$, the spectrum of $H^h$ has at least one gap in the
interval $(\alpha^h,\beta^h)$ of length larger than $Ch^{M+n}$.
\end{lem}

\begin{proof}
Divide the interval $(\alpha^h,\beta^h)$ in $[Dh^{-n}]$ equal
subintervals with some constant $D>0$ (here $[a]$ denotes the
smallest integer larger than $a$). By (\ref{l:Nh}), it follows that
if $D$ is large enough, there exists a constant $h_2>0$ such that,
for any $h\in (0,h_2]$, at least one of these intervals does not
meet the spectrum of $H^h_D$. Consider this interval and divide it
in three equal parts. By Theorem~\ref{t:D}, there exists a constant
$h_1>0$ such that, for any $h\in (0,h_1]$, the central subinterval
of this partition does not meet the spectrum of $H^h$. Since
$\alpha^h$ and $\beta^h$ belong to the spectrum of $H^h$, it is
clear that there exists a gap in $(\alpha^h,\beta^h)$, containing
this subinterval. This proves the lemma.
\end{proof}

By Lemma~\ref{l:s}, each interval $(\lambda^h_j, \lambda^h_{j+1}),
j=0,1,\ldots, N-1$ contains at least one gap in the spectrum of
$H^h$ of length $\geq Ch^{M+n}$, and the spectrum of $H^h$ on the
interval $I(h)$ has at least $N$ gaps of length $\geq Ch^{M+n}$ for
any $h$ small enough.
\end{proof}

\section{A general case}\label{s:general}
As a first application of Theorem~\ref{t:abstract}, we show that the
spectrum of the Schr\"o\-din\-ger operator $H^h$ with the periodic
magnetic field, having magnetic wells, has gaps (and, moreover, an
arbitrarily large number of gaps) on the interval $[0,
h(b_0+\epsilon_0)]$ in the semiclassical limit $h\to 0$. Under some
additional generic assumption, this result was obtained in
\cite{gaps}.

\begin{thm} \label{t:main}
Under the assumption (\ref{e:tr1}), there exists, for any natural
$N$, $h_0>0$ such that, for any $h\in (0,h_0]$, the spectrum of
$H^h$ in the interval $[0, h(b_0+\epsilon_0)]$ has at least $N$
gaps.
\end{thm}

\begin{proof}
Keep notation of Section~\ref{s:tunnel}. Fix some natural $N$.
Choose some
\[
b_0<\mu_0<\mu_1<\ldots<\mu_N< b_0+\epsilon_1.
\]
For any $j=0,1,\ldots,N$, take any $x_j\in D$ such that
\[
{\Tr}^+ (B(x_j))=\mu_j.
\]
Choose a local chart $f_j : U_j \to \RR^n$ defined in a neighborhood
$U_j$ of $x_j$ with local coordinates $X=(X_1,X_2,\ldots,X_n)\in
\RR^n$. Suppose that $f_j(U_j)$ is a ball $B=B(0,r)$ in $\RR^n$,
$f_j(x_j)=0$, the Riemannian metric at $x_j$ becomes the standard
Euclidean metric on $\RR^n$ and
\[
{\mathbf B}(x_j)=\sum_{k=1}^{d_j} \mu_{jk}dX_{2k-1}\wedge dX_{2k}.
\]
Let $\varphi_j$ be a smooth function on $B$ such that
\[
|{\mathbf A}(X)-d\varphi_j(X)-A_j^q(X)|\leq C|X|^2,
\]
where
\[
A^q_j(X)=\frac12\sum_{k=1}^{d_j}
\mu_{jk}\left(X_{2k-1}dX_{2k}-X_{2k}dX_{2k-1}\right).
\]

Write $X''=(X_{2d_j+1},\ldots,X_n)$. Let $\chi_j$ be a smooth
function on $D$ with support in a neighborhood of $x_j$ and
satisfying near $x_j$, $\chi_j(x)\equiv 1$. Let $v_j^h\in
C^\infty_c(D)$ be defined as
\[
    v_j^h(x)=\chi_j(x)\exp\left(-i\frac{\varphi_j(x)}{h}\right)
    \exp\left(-\frac{1}{4h}\sum_{k=1}^{d_j}
    \mu_{jk}(X^2_{2k-1}+X^2_{2k})\right)
    \exp\left(-\frac{|X''|^2}{h^{2/3}}\right).
\]
It is shown in the proof of Theorem 2.2 from \cite{HM} that
\[
\|(H^h_D-h\mu_j)v_j^h\|\leq Ch^{4/3}\|v_j^h\|.
\]
So the result follows from Theorem~\ref{t:abstract} with
$\mu_j^h=h\mu_j$ and $M=1$.
\end{proof}

\section{Potential wells with the regular point bottom}
One can get a more precise information on location and asymptotic
behavior of gaps in the spectrum of a magnetic Schr\"odinger
operator with magnetic wells, imposing additional conditions on the
bottom of the magnetic well. In this Section, we consider a case
when the bottom of the magnetic well contains zero-dimensional
components, that is, isolated points, and, moreover, the magnetic
field behaves regularly near these points. More precisely, we will
assume that, for some integer $k>0$, if $B(x_0)=0$, then there
exists a
  positive constant $C$ such that for all $x$ in some neighborhood
  of $x_0$ the estimate holds:
\begin{equation}\label{e:B}
C^{-1}d(x,x_0)^k\leq {\Tr}^+ (B(x))  \leq C d(x, x_0)^k
\end{equation}
(here $d(x,y)$ denotes the geodesic distance between $x$ and $y$).

In this case, the important role is played by a differential
operator $K^h_{\bar{x}_0}$ in $\RR^n$, which is in some sense an
approximation to the operator $H^h$ near $x_0$. Recall its
definition \cite{HM}.

Let $\bar x_0$ be a zero of $B$. Choose local coordinates $f: U(\bar
x_0)\to \RR^n$ on $M$, defined in a sufficiently small neighborhood
$U(\bar x_0)$ of $\bar x_0$. Suppose that $f(\bar x_0)=0$, and the
image $f(U(\bar x_0))$ is a ball $B(0,r)$ in $\RR^n$ centered at the
origin.

Write the 2-form $\bf B$ in the local coordinates as
\[
{\bf B}(X)=\sum_{1\leq l<m\leq n} b_{lm}(X)\,dX_l\wedge dX_m, \quad
X=(X_1,\ldots,X_n)\in B(0,r).
\]
Let ${\bf B}^0$ be the closed 2-form in $\RR^n$ with polynomial
components defined by the formula
\[
{\bf B}^0(X)=\sum_{1\leq l<m\leq
n}\sum_{|\alpha|=k}\frac{X^\alpha}{\alpha !}\frac{\partial^\alpha
b_{lm}}{\partial X^\alpha}(0)\,dX_l\wedge dX_m, \quad X\in\RR^n.
\]
One can find a 1-form ${\bf A}^{0}$ on $\RR^n$ with polynomial
components such that
\[
d{\bf A}^0(X) ={\bf B}^0(X), \quad X\in\RR^n.
\]

Let $K^h_{\bar{x}_0}$ be a self-adjoint differential operator in
$L^2(\RR^n)$ with polynomial coefficients given by the formula
\[
K_{\bar{x}_0}^h =  (i h\,d+{\bf A}^0)^* (i h\,d+{\bf A}^0),
\]
where the adjoints are taken with respect to the Hilbert structure
in $L^2(\RR^n)$ given by the flat Riemannian metric $(g_{lm}(0))$ in
$\RR^n$. If ${\bf A}^0$ is written as
\[
{\bf A}^0=A^0_{1}\, dX_1+\ldots+ A^0_{n}\,dX_n,
\]
then $K^h_{\bar{x}_0}$ is given by the formula
\[
K_{\bar{x}_0}^h=\sum_{1\leq l,m\leq n} g^{lm}(0) \left(i h
\frac{\partial}{\partial X_l}+A^0_{l}(X)\right)\left(i h
\frac{\partial}{\partial X_m}+A^0_{m}(X)\right).
\]
The operators $K^h_{\bar{x}_0}$ have discrete spectrum (cf, for
instance, \cite{HelNo85,HM88}). Using the simple dilation $X\mapsto
h^{\frac{1}{k+2}}X$, one can show that the operator
$K^h_{\bar{x}_0}$ is unitarily equivalent to
$h^{\frac{2k+2}{k+2}}K^1_{\bar{x}_0}$. Thus,
$h^{-\frac{2k+2}{k+2}}K^h_{\bar{x}_0}$ has discrete spectrum,
independent of $h$.

\begin{thm}\label{t:gaps0}
Suppose that the operator $H^h$ satisfies the condition
(\ref{e:tr1}) with some $\epsilon_0>0$ and there exists a zero
$\bar{x}_0$ of $B$, $B(\bar{x}_0)=0$, satisfying the assumption
(\ref{e:B}) for some integer $k>0$. Then, for any natural $N$, there
exist constants $C>0$ and $h_0>0$ such that the part of the spectrum
of $H^h$, contained in the interval $[0,Ch^{\frac{2k+2}{k+2}}]$, has
at least $N$ gaps for any $h\in (0,h_0)$.
\end{thm}

\begin{proof}
Fix $\epsilon_1$ and $\epsilon_2$ such that $0<\epsilon_1 <
\epsilon_2 < \epsilon_0$ and consider the operator $H^h_D$
associated with the domain $D=\overline{U_{\epsilon_2}}$. Denote by
$\lambda_1<\lambda_2<\lambda_2<\ldots$ the eigenvalues of the
operator $K^1_{\bar{x}_0}$ (not taking into account multiplicities).

For any $j\in \NN$, let $w^h_j\in L^2(\RR^n)$ be any eigenfunction
of $K^h_{\bar{x}_0}$ corresponding to the eigenvalue
$h^{\frac{2k+2}{k+2}}\lambda_j$. Let $\chi$ be a compactly supported
cut-off function in the neighborhood $U(\bar x_0)$ of $\bar{x}_0$ as
above equal to $1$ in a neighborhood of $\bar{x}_0$. Define
\[
v^h_j(x)=\chi(x)w^h_j(x).
\]
As shown in the proof of Theorem 2.5 in \cite{HM}, we have
\[
\|\left(H^h_D-h^{\frac{2k+2}{k+2}}\lambda_j\right)v^h_j\|\leq C_j
h^{\frac{2k+3}{k+2}}\|v^h_j\|.
\]
For a given natural $N$, choose any constant $C>\lambda_{N+1}$. Then
the result follows from Theorem~\ref{t:abstract} with
$\mu_j^h=h^{\frac{2k+2}{k+2}}\lambda_j, j=1,\ldots,N+1$.
\end{proof}

\section{Potential wells with the one-dimensional bottom}
In this section we consider the case when the manifold $M$ is an
oriented two-dimensional Riemannian manifold, and the zero set of
the periodic magnetic field $\mathbf B$ contains a one-dimensional
non-degenerate compact component. More precisely, suppose that:
\begin{itemize}
  \item $b_0=0$;
  \item the zero set of the magnetic field $\mathbf B$ has a connected component
  $\gamma$, which is a bounded smooth curve;
\item there are constants $k\in \NN$ and $C>0$ such that for all $x$
in some neighborhood of $\gamma$ the estimates hold:
\begin{equation}\label{e:B1}
C^{-1}d(x,\gamma)^k\leq |B(x)|  \leq C d(x,\gamma)^k\,.
\end{equation}
\end{itemize}
In particular, for $k=1$ the last condition means that $\nabla B$
does not vanish on $\gamma$.

On compact manifolds, this model was introduced for the first time
by Montgomery \cite{Mont} and was further studied in
\cite{HM,Pan-Kwek,syrievienne}. In this paper we consider the
problem of existence of gaps for this model in a simple particular
case. Namely, we assume that the leading term of the Taylor
expansion of the magnetic field $B$ at $\gamma$ is constant along
$\gamma$. More precisely, write the volume 2-form ${\mathbf B}$ as
\[
{\mathbf B}=B(x)\omega, \quad x\in M,
\]
where $\omega$ is the Riemannian volume form on $M$. Denote by $N$
the external unit normal vector to $\gamma$. Let $\tilde{N}$ denote
an arbitrary extension of $N$ to a smooth vector field on $M$.
Consider the function $\beta_1$ on $M$ given by the formula
\[
\beta_1(x)=\tilde{N}^kB(x), \quad x\in M.
\]
By (\ref{e:B1}), it is easy to see that
\[
\beta_1(x)\not=0, \quad x\in \gamma.
\]
Our assumption is that the restriction of $\beta_1$ to $\gamma$
(which is independent of the choice of smooth extension
$\tilde{N}$) is constant along $\gamma$:
\begin{equation}\label{e:B2}
\beta_1(x)\equiv \beta_1 = {\rm const}, \quad x\in \gamma.
\end{equation}
For $k=1$, this condition means that the length of the gradient
$|\nabla B|$ is constant along $\gamma$.

\begin{thm}
Under the given assumptions, for any natural $N$ there exist
constants $C>0$ and $h_0>0$ such that the part of the spectrum of
$H^h$ contained in the interval $[0,Ch^{\frac{2k+2}{k+2}}]$ has at
least $N$ gaps for any $h\in (0,h_0)$.
\end{thm}

\begin{proof}
Choose a normal coordinate system $(x,y)$ in a tubular
neighborhood $U$ of $\gamma$. Thus, $x\in [0,L)\cong
S^1_L=\RR/L\ZZ$ is the natural parameter along $\gamma$ ($L$ is
the length of $\gamma$), $\gamma$ is given by the equation $y=0$,
and $y\in (-\varepsilon_0,\varepsilon_0)$ is the natural parameter
along the geodesic, passing through the point on $\gamma$ with the
coordinates $(x,0)$ orthogonal to $\gamma$. It is well known that
in such coordinates the metric $g$ has the form
\[
g=a(x,y)^2dx^2+dy^2,
\]
where
\[
a(x,0)=1,\quad \frac{\partial a}{\partial y}(x,0)=0.
\]
Write
\[
{\bf A}=A_1dx+A_2dy
\]
and
\[
{\bf B}=b(x,y)dx\wedge dy, \quad b=\frac{\partial A_2}{\partial
x}-\frac{\partial A_1}{\partial y}.
\]
The external unit normal vector to $\gamma$ has the form
\[
N=\frac{\partial }{\partial y},
\]
and one can take a smooth extension $\tilde{N}$ as
\[
\tilde{N}=\frac{\partial }{\partial y}.
\]
Thus, we have
\[
\frac{\partial^j b}{\partial y^j}(x,0)=0, \quad j=0,1,\ldots, k-1.
\]
and
\[
\beta_1=\left.\frac{\partial^k}{\partial y^k }
\left(a^{-1}b\right)\right|_{y=0}=\frac{\partial^k b}{\partial
y^k}(x,0)\not=0.
\]

As above, denote by $H_D^h$ the unbounded self-adjoint operator in
$L^2(D)$ given by the operator $H^h$ in the domain $D=\overline{U}$
with the Dirichlet boundary conditions and by $\lambda_1^h<
\lambda_2^h< \ldots < \lambda^h_{N}<\ldots$ the eigenvalues of
$H_D^h$.

Adding an exact one form $d\phi(x)$ to ${\bf A}$, without loss of
generality, we can assume that
\[
A_1(x,0)=\alpha_1 \equiv {\rm const}.
\]

Consider the self-adjoint operator $H^{h,0}$ in $L^2(S^1_L\times
\RR)$ defined by the formula
\[
H^{h,0}=-h^2\frac{\partial^2 }{\partial y^2}+ \left(ih\frac{\partial
}{\partial x}+\alpha_1+\frac{1}{(k+1)!}\beta_1 y^{k+1}\right)^2,
\quad x\in S^1_L, \quad y\in \RR.
\]

By \cite[Theorem 2.7]{HM}, the operator $H^{h,0}$ in $L^2 (S^1_L
\times \RR)$ has discrete spectrum. We construct some eigenfunctions
of $H^{h,0}$, using separation of variables. Consider a function
$u\in L^2(S^1_L\times \RR)$ of the form
\[
u(x,y)=e^{2\pi i\frac{p(h)}{L}x}v(y), \quad x\in S^1_L, \quad y\in
\RR,
\]
with some $v\in L^2(\RR, dy)\cap C^\infty(\RR)$ and $p(h)\in \ZZ$.
Then
\[
H^{h,0}u(x,y)=e^{2\pi i\frac{p(h)}{L}x}H\left(h,\beta(h)\right)v(y),
\]
where
\[
\beta(h)=\frac{2\pi h p(h)}{L}-\alpha_1
\]
and
\[
H(h,\beta)=-h^2\frac{\partial^2 }{\partial y^2}+
\left(\beta-\frac{1}{(k+1)!}y^{k+1}\right)^2.
\]
For any $\alpha>0$ the dilation operator
\[
(T(\alpha)f)(y)=\sqrt{\alpha}f(\alpha y), \quad f\in L^2(\RR, dy),
\]
is a unitary operator in $L^2(\RR,dy)$, satisfying the conditions
\[
\frac{\partial}{\partial y}T(\alpha)=\alpha
T(\alpha)\frac{\partial}{\partial y}, \quad
y^sT(\alpha)=\alpha^{-s} T(\alpha)y^s.
\]
Using these relations, it is easy to check that the identity
\[
H(h,\beta)T(\alpha)=\alpha^{-(2k+2)}T(\alpha)
H(\alpha^{k+2}h,\alpha^{k+1}\beta)
\]
holds for any $h>0, \alpha>0, \beta>0$, and, in particular,
\[
H(h,\beta)T(h^{-\frac{1}{k+2}})=h^{\frac{2k+2}{k+2}}T(h^{-\frac{1}{k+2}})
H(1,h^{-\frac{k+1}{k+2}}\beta).
\]
For any fixed $b\in \RR$, the operator $H(1,b)$ has simple
discrete spectrum
\[
\mu_1(b) < \mu_2(b)<\ldots, \quad \mu_j(b)\to +\infty,
\]
where $\mu_j(b)$ are continuous functions. Fix a natural $N$. Then
there exist an interval $(b_1,b_2)$ and a system of disjoint
intervals $(c_j,C_j), j=1,2,\ldots,N$, $c_1<C_1<c_2<\ldots <
C_{N-1}<c_N<C_N$ such that, for any $b\in (b_1,b_2)$, the
inclusions $\mu_j(b)\in (c_j,C_j), j=1,2,\ldots,N$ hold. Choose
$p(h)\in \ZZ$ so that
\[
b_1<h^{-\frac{k+1}{k+2}} \beta(h) <b_2,
\]
or, equivalently,
\[
\frac{L}{2\pi} (\alpha_1h^{-1}+b_1h^{-\frac{1}{k+2}}) <p(h)
<\frac{L}{2\pi} (\alpha_1h^{-1}+b_2h^{-\frac{1}{k+2}}).
\]
Such a $p(h)$, clearly, always exists.

Let $v_j\in L^2(\RR, dy)\cap C^\infty(\RR)$ be the eigenfunction
of $H(1,h^{-1/(k+2)}\beta(h))$, corresponding to the eigenvalue
$\mu_j(h^{-1/(k+2)}\beta(h))$:
\[
H\left(1,h^{-1/(k+2)}\beta(h)\right)v_j=
\mu_j\left(h^{-1/(k+2)}\beta(h)\right)v_j.
\]
Then $v_j^h=T(h^{-1/(k+2)})v_j$ is an eigenfunction of
$H\left(h,\beta(h)\right)$ with the corresponding eigenvalue
\[
\mu^h_j=h^{\frac{2k+2}{k+2}}\mu_j\left(\beta(h)\right),
\]
and, therefore, the function
\[
u^h_j(x,y)=e^{2\pi i\frac{p(h)}{L}x}v^h_j(y)
\]
is an eigenfunction of $H^{h,0}$ with the same eigenvalue:
\[
H^{h,0}u^h_j=\mu^h_j u^h_j.
\]
Thus, we have shown that, for any natural $N$, there exists $N$
eigenvalues $\mu^h_{1}< \mu^h_{2}< \ldots < \mu^h_{N}$ of $H^{h,0}$
(generally speaking, not consecutive) such that
\[
\mu^h_{j+1}-\mu^h_{j}>Ch^{\frac{2k+2}{k+2}}, \quad j=1,2,\ldots,N-1.
\]

The normal coordinate system $(x,y)$ defined above gives a
diffeomorphism
\[
\Theta:  S^1_L \times (-\varepsilon_0,\varepsilon_0) \to \Omega,
\]
onto a tubular neighborhood $\Omega$ of $\gamma$. Let $\chi_0$ be a
cut off function supported in $S^1_L \times
(-\varepsilon_0,\varepsilon_0)$ and equal to $1$ in a neighborhood
of $S^1_L\times \{0\}$. Then, by \cite[Theorem 2.7]{HM}, we have
\[
\|(H_D^h-\mu^h_j) (\chi_0 u^h_j)\circ \Theta^{-1}\|\leq C
h^{\frac{2k+2}{k+2}}\|\chi_0 u^h_j\|, \quad j=1,2,\ldots,N.
\]
Theorem~\ref{t:abstract} completes the proof.
\end{proof}

\section{Concluding remarks}
In \cite{Ko04}, one considered the case when all the zeroes of $B$
are isolated points, satisfying the assumption (\ref{e:B}) for some
integer $k>0$. In this case one can get a more precise information
about location of spectral gaps.

\begin{thm}[\cite{Ko04}]\label{t:cmp}
Suppose that there exists at least one zero of $B$, and all the
zeroes of $B$ are isolated points, satisfying the assumption
(\ref{e:B}) for some integer $k>0$. Then there exists an increasing
sequence $\{\lambda_m, m\in\NN \}$, satisfying the condition:
$\lambda_m\to\infty$ as $m\to\infty$, such that, for any $a$ and $b$
such that $\lambda_m<a<b < \lambda_{m+1}$ for some $m$,
\[
[ah^{\frac{2k+2}{k+2}}, bh^{\frac{2k+2}{k+2}}]\cap
\sigma(H^h)=\emptyset\,,
\]
for any $h>0$ small enough.
\end{thm}

As a direct consequence of Theorem~\ref{t:cmp}, we get another proof
of Theorem~\ref{t:gaps0} in this particular case.

The numbers $\{\lambda_m, m\in\NN \}$ in Theorem~\ref{t:cmp} are the
eigenvalues of a so called model operator $K^h$, which is defined as
follows. Choose a fundamental domain $\cF\subset M$ so that $B$ does
not vanishes on the boundary of $\cF$. Let $\{\bar
x_j|\,j=1,\dots,K\}$ denote all the zeroes of $B$ in $\cF$; $\bar
x_i\ne\bar x_j$, if $i\ne j$. The operator $K^h$ is a self-adjoint
operator in $L^2({\mathbb R}^n)^K$ given by
\[
K^h = \bigoplus_{1\le j\le K} K^h_{\bar{x}_j}.
\]

The proof of Theorem~\ref{t:cmp} given in \cite{Ko04} makes use of
abstract functional-analytic methods developed in~\cite{KMS} (see
also a survey paper \cite{bedlewo}). These methods were developed
for the study of similar questions for a periodic magnetic
Schr\"odinger operator
\[
H_\mu=(i\,d+{\bf A})^* (i\,d+{\bf A})+\mu^{-2}V(x),
\]
with a $\Gamma$-invariant Morse potential $V\geq 0$ on the universal
covering $M$ of a compact manifold in the strong electric field
limit ($\mu\to 0$), where, as above, we assume that the 2-form ${\bf
B}=d{\bf A}$ is $\Gamma$-invariant (here $\Gamma=\pi_1(M)$). Indeed,
they allow to obtain stronger results than the existence of gaps in
the spectrum, namely, to prove Murrey-von Neumann equivalence of the
corresponding spectral projections of $H_\mu$ and the associated
model operator (see \cite{KMS,bedlewo} for details).

Observe, however, that, using Theorem~\ref{t:D} and making full use
of Theorem 2.5 in \cite{HM} (that provides a more complete
information on the whole spectrum of $H^h_D$, not only the
construction of some approximate eigenvalues as it was needed for
the proof of Theorem~\ref{t:gaps0}), one can easily give another
proof of Theorem~\ref{t:cmp}.

On the other side, it is much easier to construct approximate
eigenvalues, showing the existence of some eigenvalues for $H^h_D$,
than to localize completely its spectrum. Therefore, the scheme of a
proof of existence of spectral gaps suggested in this paper could be
very efficient for treating cases where we know how to construct
approximate eigenvalues, but the treatment of the whole spectrum
could be more difficult. One could mention two cases where this idea
should definitely work:
\begin{enumerate}
  \item The two-dimensional magnetic Dirichlet problem studied in \cite{HM01}
  (see, in particular, \cite[Theorem 7.3]{HM01}).
  \item The more generic Montgomery model considered in \cite{syrievienne}.
\end{enumerate}
These cases will be discussed elsewhere.


\begin{thebibliography}{99}
\bibitem{Ag} S. Agmon, {\em Lectures on exponential decay of
solutions of second-order elliptic equations.} Mathematical Notes
29, Princeton University Press, Princeton, 1982.

\bibitem{BDP} J. Br\"uning, S. Yu. Dobrokhotov, K. V. Pankrashkin. The
spectral asymptotics of the two-dimensional Schr\"odinger operator
with a strong magnetic field. I. Russ. J. Math. Phys. 9 (2002), no.
1, 14--49; II. Russ. J. Math. Phys. 9 (2002), no. 4, 400--416 (see
also e-print version math-ph/0411012).

\bibitem{Carlsson} U. Carlsson, An infinite number of wells in the
semi-classical limit. {\em Asymptotic Anal.} {\bf 3} (1990), no.
3, 189--214.

\bibitem{Di-Sj} M. Dimassi, J. Sj\"ostrand, {\em Spectral
asymptotics in the semi-classical limit.} London Mathematical
Society Lecture Notes Series, 268, Cambridge University Press,
Cam\-bri\-dge, 1999.

\bibitem{Frank} R. L. Frank, On the tunneling effect for magnetic
Schr\"odinger operators in antidot lattices. {\em Asymptot. Anal.}
{\bf 48} (2006), no. 1-2, 91--120.

\bibitem{syrievienne} B. Helffer, Introduction to semi-classical methods
for the Schr\"odinger operator with magnetic fields. To appear in:
S\'eminaires et Congr\`es, Cours du CIMPA, SMF.

\bibitem{gaps} B. Helffer, Yu. A. Kordyukov, Semiclassical
asymptotics and gaps in the spectra of periodic Schr\"odinger
operators with magnetic wells, preprint math.SP/0601366; to appear
in {\em Trans. Amer. Math. Soc.}

\bibitem{HM88} B. Helffer, A. Mohamed, Caract\'erisation du
spectre essentiel de l'op\'erateur de Schr\"odinger avec un champ
magn\'etique. {\em Ann. Inst. Fourier} {\bf 38} (1988) 95--112.

\bibitem{HM} B. Helffer, A. Mohamed, Semiclassical analysis
for the ground state energy of a Schr\"odinger operator with
magnetic wells. {\em J. Funct. Anal.} {\bf 138} (1996), 40--81.

\bibitem{HM01} B. Helffer, A. Morame, Magnetic bottles in connection with
superconductivity. {\em J. Funct. Anal.} {\bf 185} (2001), 604--680.

\bibitem{HelNo85} B. Helffer, J. Nourrigat, {\em
Hypoellipticit\'e maximale pour des op\'erateurs polyn\^{o}mes de
champs de vecteurs.} Boston: Birkh\"auser, 1985

\bibitem{HSI}
B. Helffer, J. Sj\"ostrand,  Multiple wells in the semiclassical
limit. I. {\em Comm. Partial Differential Equations}  \textbf{9}
(1984),  337--408.

\bibitem{HSII}
B. Helffer, J. Sj\"ostrand,  Puits multiples en limite
semi-classique. II. Interaction mol\'eculaire. Sym\'etries.
Perturbation. {\em Ann. Inst. H. Poincar\'e Phys. Th\'eor.}
\textbf{42}, no. 2 (1985), 127--212.

\bibitem{HS87} B. Helffer, J. Sj\"ostrand,
Effet tunnel pour l'\'equation de Schr\"odinger avec champ
magn\'etique. {\em Ann. Scuola Norm. Sup. Pisa Cl. Sci.} (4)
\textbf{14}  (1987),  625--657.

\bibitem{HS88} B. Helffer, J. Sj\"ostrand,
Analyse semi-classique pour l'\'equation de Harper (avec
application \`a  l'\'equation de Schr\"odinger avec champ
magn\'etique). {\em M\'em. Soc. Math. France (N.S.)}  \textbf{34}
(1988).


\bibitem{HSLNP345} B. Helffer, J. Sj\"ostrand, \'{E}quation de
{S}chr\"odinger avec champ magn\'etique et \'equation de {H}arper,
In: {\em Schr\"odinger operators (S\o nderborg, 1988)}, Lecture
Notes in Phys., 345, Springer, Berlin, 1989, pp. 118--197.

\bibitem{HempelHerbst95} R. Hempel, I. Herbst, Strong magnetic
fields, Dirichlet boundaries, and spectral gaps. {\em Commun. Math.
Phys.} {\bf 169} (1995), 237--259.

\bibitem{HempelPost02} R. Hempel, O. Post, Spectral gaps for
periodic elliptic operators with high contrast: an overview, In:
{\em Progress in analysis, Vol. I, II (Berlin, 2001).} World Sci.
Publishing, River Edge, NJ, 2003, pp. 577--587.

\bibitem{HerbstNakamura}
I. Herbst, S. Nakamura. Schrodinger operators with strong magnetic
fields: quasi-periodicity of spectral orbits and topology. In: {\em
Differential operators and spectral theory, Amer. Math. Soc. Transl.
Ser. 2, 189}, Amer. Math. Soc., Providence, RI, 1999, pp. 105--123.

\bibitem{KMS} Yu. A. Kordyukov, V. Mathai, M. Shubin,
Equivalence of projections  in semiclassical limit and a vanishing
theorem for higher traces in $K$-theory. {\em J. Reine Angew.
Math.} {\bf 581} (2005), 193--236.

\bibitem{Ko04} Yu. A. Kordyukov,
Spectral gaps for periodic Schr\"odinger operators with strong
magnetic fields. {\em Commun. Math. Phys.} {\bf 253} (2005),
371--384.

\bibitem{bedlewo} Yu. A. Kordyukov, Semiclassical asymptotics and spectral gaps for periodic
magnetic Schr\"odinger operators on covering manifolds, In: {\em
``$C^*$-algebras and elliptic theory'', Trends in Mathematics},
129 -- 150, Birkh\"auser, Basel, 2006.


\bibitem{MS} V. Mathai, M. Shubin, Semiclassical asymptotics
and gaps in the spectra of magnetic Schr\"odinger operators. {\em
Geometriae Dedicata} {\bf 91} (2002), 155--173.

\bibitem{Mont} R. Montgomery, Hearing the zero locus of a magnetic field. {\em
Comm. Math. Phys.} {\bf 168} (1995), 651--675.


\bibitem{Nakamura95} S. Nakamura, Band spectrum for Schr\"odinger
operators with strong periodic magnetic fields. In: {\em Partial
differential operators and mathematical physics (Holzhau, 1994),
Operator Theory: Advances and Applications. vol. 78}, 261--270,
Birkh\"auser, Basel, 1995.

\bibitem{Pan-Kwek} X-B. Pan, K-H. Kwek. Schr\"odinger operators with
non-degenerately vanishing magnetic fields in bounded domains.
Trans. Amer. Math. Soc. {\bf 354} (2002), 4201--4227.
\end{thebibliography}
\end{document}